\documentclass{article}
\usepackage{amssymb}

\usepackage{graphicx}
\usepackage{amsmath}

\newtheorem{theorem}{Theorem}

\newtheorem{corollary}[theorem]{Corollary}

\newtheorem{lemma}[theorem]{Lemma}

\newtheorem{proposition}[theorem]{Proposition}

\newenvironment{proof}[1][Proof]{\textbf{#1.} }{\ \rule{0.5em}{0.5em}}
\input{tcilatex}

\begin{document}

\title{Another Coboundary Operator for Differential Forms with Values in the Lie
Algebra Bundle of a Group Bundle\\
$\lhd$A Chapter in Synthetic Differential Geometry of Groupoids$\rhd$}
\author{Hirokazu Nishimura \\
Institute of Mathematics, University of Tsukuba\\
Tsukuba, Ibaraki, 305-8571, Japan}
\maketitle

\begin{abstract}
Kock [Bull. Austral. Math. Soc., 25 (1982), 357-386] has considered
differential forms with values in a group in a context where neighborhood
relations are available. By doing so, he has made it clear where the
so-called Maurer-Cartan formula should come from. In this paper, while we
retain the classical definition of differential form with values in the Lie
algebra of a group, we propose another definition of coboundary operator for
the de Rham complex in a highly general microlinear context, in which
neighborhood relations are no longer in view. Using this new definition of
coboundary operator, it is to be shown that the main result of Kock's paper
mentioned above still prevails in our general microlinear context. Our
considerations will be carried out within the framework of groupoids.
\end{abstract}

\section{\label{s0}Introduction}

Although the Maurer-Cartan equation has long been known, it is Kock \cite{k3}
that excavated its hidden geometric meaning for the first time. To this end,
he introduced differential forms with values in groups in place of classical
ones with values in their Lie algebras by using first neighborhood
relations. We do not know exactly what spaces enjoy first neighborhood
relations, but we are sure at least that formal manifolds are naturally
endowed with such relations. As far as formal manifolds are concerned, he
has demonstrated that his new definition of differential form is anyway
equivalent to the classical one. This means that if we want to generalize
his geometric ideas from formal manifolds to microlinear spaces in general,
it is not necessary to adhere to his noble definition of differential form.
In other words, we can say that his novel definition of differential form
does not consitute the indispensible components of his thrilling geometic
ideas.

Another unsatisfactory feature of \cite{k3} is that his proof on the exact
discrepancy between the coboundary operator of the de Rham complex and his
contour derivative from dimension $1$ to dimension $2$, from which the
Maurer-Cartan equation comes at once, appears rather analytic than
synthetic, though he is very famous in synthetic differential geometry.

We are inexpiably diehard in the definition of differential form, so that we
prefer differential forms with values in their Lie algebras to ones with
values in groups themselves. We prefer cubical arguments to simplicial ones
which Kock admires. However we propose another definition for the coboundary
operator of the de Rham complex, which is to be eventually shown to be
equivalent to the classical one in synthetic differential geometry. This new
definition of the coboundary operator facilitates the comparision between
the coboundary operator of the de Rham complex and the cubical version of
his contour derivative from dimension $1$ to dimension $2$. Thus this paper
might be put down as a microlinear generalization of Kock's \cite{k3}
ingeneous ideas.

The paper is organized as follows. After giving some preliminaries and
fixing our notation in the succeeding section, we will deal with the
coboundary operator for the de Rham complex in a somewhat general context,
in which differential forms on a groupoid $G$ over a base space $M$ with
values in a vector bundle $E$ over the same base space $M$ are considered.
This is the topic of \S \ref{s2}. The cobondary operator studied here is
called the \textit{additive coboundary operator} and is denoted $\mathbf{d}%
_{+}$. If the vector bundle $E$ happens to be the Lie algebra bundle $%
\mathcal{A}^{1}L$ of a group bundle $L$ over $M$, another definition of the
coboundary operator naturally emerges besides the additive one with due
regard to group structures that $L$ possesses. The emerging coboundary
operator is called the \textit{multiplicative coboundary operator} and is
denoted $\mathbf{d}_{\times }$. This is the topic of \S \ref{s3}. The
succeeding section is devoted to establishing the coincidence of $\mathbf{d}%
_{+}$ and $\mathbf{d}_{\times }$ whenever both are available. The last
section is devoted to a microlinear generalization of Kock's \cite{k3} main
result. Our standard reference on synthetic differential geometry is
Lavendhomme \cite{l1}. Unless stated to the contary, every space in this
paper is assumed to be microlinear. Our discussions will be carried out
within the context of groupoids, which is a bit more general than Kock's 
\cite{k3}.

\section{\label{s1}Preliminaries}

\subsection{\label{p1}Groupoids}

Let $M$ be a microlinear space. Given $x\in M$ and a groupoid $G$ over a
base $M$ with its object inclusion map $\mathrm{id}:M\rightarrow G$ and its
source and target projections $\alpha,\beta:G\rightarrow M$, we denote by $%
\mathcal{A}_{x}^{n}G$ the totality of mappings $\gamma:D^{n}\rightarrow G$
with $\gamma(0,...,0)=\mathrm{id}_{x}$ and $(\alpha\circ%
\gamma)(d_{1},...,d_{n})=x$ for any $(d_{1},...,d_{n})\in D^{n}$. We denote
by $\mathcal{A}^{n}G$ the set-theoretic union of $\mathcal{A}_{x}^{n}G$'s
for all $x\in M$. Given $\gamma\in\mathcal{A}^{n}G$ and $(d_{1},...,d_{n})%
\in D^{n}$, we will often write $\gamma_{d_{1},...,d_{n}}$ in place of $%
\gamma (d_{1},...,d_{n})$. If $G$ is the pair groupoid (i.e., $G=M\times M$%
), then $\mathcal{A}^{n}G$ can and should be identified with $M^{D^{n}}$.
The canonical projection $\pi_{n}:\mathcal{A}^{n}G\rightarrow M$ is defined
to be 
\begin{equation*}
\pi_{n}(\gamma)=\gamma_{0,...,0}
\end{equation*}
for any $\gamma\in\mathcal{A}^{n}G$. The canonical mapping $\mathbf{s}_{i}:%
\mathcal{A}^{n}G\rightarrow\mathcal{A}^{n+1}G$ ($i=1,...,n+1$) is defined to
be 
\begin{equation*}
\mathbf{s}_{i}(\gamma)_{d_{1},...,d_{n+1}}=\gamma_{d_{1},...,\widehat{d_{i}}%
,...,d_{n+1}}
\end{equation*}
for any $\gamma\in\mathcal{A}^{n}G$ and any $(d_{1},...d_{n+1})\in D^{n+1}$.
The canonical mapping $\mathbf{d}_{i}:\mathcal{A}^{n+1}G\rightarrow \mathcal{%
A}^{n}G$ ($i=1,...,n+1$) is defined to be 
\begin{equation*}
\mathbf{d}_{i}(\gamma)_{d_{1},...,d_{n}}=%
\gamma_{d_{1},...,d_{i-1},0,d_{i},...,d_{n}}
\end{equation*}
for any $\gamma\in\mathcal{A}^{n+1}G$ and any $(d_{1},...,d_{n})\in D^{n}$.
It is easy to see that $\mathbf{s}_{i}$'s and $\mathbf{d}_{i}$'s satisfy the
so-called simplicial identities. By the same token as in Propositions 1 and
2 of \S3.1 in \cite{l1}, it is easy to see that $\mathcal{A}^{1}G$ can
naturally be regarded as a vector bundle over $M$ (i.e., $\mathcal{A}%
_{x}^{1}G$ is a Euclidean $\mathbb{R}$-module for any $x\in M$, where $%
\mathbb{R}$ stands for the set of real numbers with a cornucopia of
nilpotent infinitesimals pursuant to the general Kock-Lawvere axiom).
Similarly the square 
\begin{equation*}
\begin{array}{ccc}
\mathcal{A}^{2}G & 
\begin{array}{c}
\mathbf{d}_{1} \\ 
\rightarrow
\end{array}
& \mathcal{A}^{1}G \\ 
\begin{array}{cc}
\mathbf{d}_{2} & \downarrow
\end{array}
&  & 
\begin{array}{cc}
\downarrow & \pi_{1}
\end{array}
\\ 
\mathcal{A}^{1}G & 
\begin{array}{c}
\mathbf{\rightarrow} \\ 
\pi_{1}
\end{array}
& M
\end{array}
\end{equation*}
is a double vector bundle over $M$ in the sense of Mackenzie \cite{ma1},
\S9.1. The canonical mapping $\mathbf{a}_{n}:\mathcal{A}^{n}G\rightarrow
M^{D^{n}}$ is defined to be 
\begin{equation*}
\mathbf{a}_{n}(\gamma)_{d_{1},...,d_{n}}=\beta(\gamma_{d_{1},...,d_{n}}) 
\end{equation*}
for any $\gamma\in\mathcal{A}^{n}G$ and any $(d_{1},...,d_{n})\in D^{n}$.

Given a group bundle $L$ over $M$, $\mathcal{A}L$ is not only a vector
bundle over $M$ but also a Lie algebra bundle over $M$, where the Lie
bracket $[\cdot,\cdot]$ is defined by the following proposition.

\begin{proposition}
\label{t1.1.1}Let $x\in M$. Given $t_{1},t_{2}\in\mathcal{A}_{x}^{1}L$,
there exists a unique $[t_{1},t_{2}]\in\mathcal{A}_{x}^{1}L$ such that 
\begin{equation}
\lbrack
t_{1},t_{2}]_{d_{1}d_{2}}=(t_{2})_{-d_{2}}(t_{1})_{-d_{1}}(t_{2})_{d_{2}}(t_{1})_{d_{1}} 
\label{1.1.1}
\end{equation}
for any $d_{1},d_{2}\in D$.
\end{proposition}

\begin{proof}
By Proposition 7 of \S2.2 in \cite{l1}, it suffices to note that if $d_{1}=0 
$ or $d_{2}=0$, then the right-hand side of (\ref{1.1.1}) is equal to $%
\mathrm{id}_{x}$, which is easy to see.
\end{proof}

\begin{theorem}
\label{t1.1.2}With respect to $[\cdot,\cdot]$ defined above, $\mathcal{A}%
_{x}^{1}L$ is a Lie algebra.
\end{theorem}

\begin{proof}
By the same token as in our \cite{n2}.
\end{proof}

We note the following simple proposition in passing.

\begin{proposition}
\label{t1.1.3}Let $x\in M$. Given $t_{1},t_{2}\in\mathcal{A}_{x}^{1}L$, we
have 
\begin{equation*}
(t_{1}+t_{2})_{d}=(t_{2})_{d}(t_{1})_{d}=(t_{1})_{d}(t_{2})_{d}
\end{equation*}
for any $d\in D$.
\end{proposition}

\begin{proof}
By the same token as in Proposition 7 of \cite{n2}.
\end{proof}

As an easy corollary of this proposition, we can see that 
\begin{equation*}
t_{-d}=(t_{d})^{-1}
\end{equation*}
for any $t\in\mathcal{A}_{x}^{1}L$ and any $d\in D$, since we have $%
(d,-d)\in D(2)$.

Our standard reference on groupoids is \cite{ma1}.

\subsection{Differential Forms}

Given a groupoid $G$ and a vector bundle $E$ over the same space $M$, the
space $\mathbf{C}^{n}(G,E)$ of \textit{differential }$n$\textit{-forms with
values in} $E$ consists of all mappings $\omega$ from $\mathcal{A}^{n}G$ to $%
E$ whose restriction to $\mathcal{A}_{x}^{n}G$ for each $x\in M$ takes
values in $E_{x}$ satisfying the following $n$-homogeneous and alternating
properties:

\begin{enumerate}
\item  We have 
\begin{equation*}
\omega(a\underset{i}{\cdot}\gamma)=a\omega(\gamma)\text{ \ \ \ \ \ }(1\leq
i\leq n) 
\end{equation*}
for any $a\in\mathbb{R}$ and any $\gamma\in\mathcal{A}_{x}^{n}G$, where $a%
\underset{i}{\cdot}\gamma\in\mathcal{A}_{x}^{n}G$ is defined to be 
\begin{equation*}
(a\underset{i}{\cdot}\gamma)_{d_{1},...,d_{n}}=%
\gamma_{d_{1},...,d_{i-1},ad_{i},d_{i+1},...d_{n}}
\end{equation*}
for any $(d_{1},...,d_{n})\in D^{n}$.

\item  We have 
\begin{equation*}
\omega(\gamma\circ D^{\sigma})=\mathrm{sign}(\sigma)\omega(\gamma) 
\end{equation*}
for any permutation $\sigma$ of $\{1,...,n\}$, where $D^{\sigma}:D^{n}%
\rightarrow D^{n}$ permutes the $n$ coordinates by $\sigma$.
\end{enumerate}

\subsection{Two Infinitesimal Stokes' Theorems}

Let $E$ be a vector bundle over $M$. If $\omega\in\mathbf{C}^{n}(G,E)$, then
the mapping $\varphi_{\omega}:\mathcal{A}^{n}G\times D^{n}\rightarrow E$
defined by 
\begin{equation*}
\varphi_{\omega}(\gamma,d_{1},...,d_{n})=d_{1}...d_{n}\omega(\gamma) 
\end{equation*}
abides by the following conditions:

\begin{enumerate}
\item  We have 
\begin{equation*}
\varphi_{\omega}(a\underset{i}{\cdot}\gamma,d_{1},...,d_{n})=a\varphi_{%
\omega }(\gamma,d_{1},...,d_{n})\text{ \ \ \ \ }(1\leq i\leq n) 
\end{equation*}
for any $a\in\mathbb{R}$.

\item  We have 
\begin{align*}
& \varphi_{\omega}(\gamma,d_{1},...,d_{i-1},ad_{i},d_{i+1},...d_{n}) \\
& =a\varphi_{\omega}(\gamma,d_{1},...,d_{i-1},d_{i},d_{i+1},...d_{n})\text{
\ \ \ \ }(1\leq i\leq n)
\end{align*}
for any $a\in\mathbb{R}$.

\item  We have 
\begin{equation*}
\varphi_{\omega}(\gamma\circ D^{\sigma},d_{1},...,d_{n})=\mathrm{sign}%
(\sigma)\varphi_{\omega}(\gamma,d_{1},...,d_{n}) 
\end{equation*}
for any permutation $\sigma$ of $\{1,...,n\}$.
\end{enumerate}

Conversely we have

\begin{theorem}
\label{t1.3.1}If $\varphi:\mathcal{A}^{n}G\times D^{n}\rightarrow E$
satisfies the above three conditions, then there exists a unique $%
\omega_{\varphi}\in\mathbf{C}^{n}(G,E)$ such that 
\begin{equation*}
\varphi(\gamma,d_{1},...,d_{n})=d_{1}...d_{n}\omega_{\varphi}(\gamma) 
\end{equation*}
for any $\gamma\in\mathcal{A}^{n}G$ and any $(d_{1},...,d_{n})\in D^{n}$.
\end{theorem}

\begin{proof}
By the same token as in the proof of Proposition 2 of \S4.2 of Lavendhomme 
\cite{l1}.
\end{proof}

Let $L$ be a group bundle over $M$. If $\omega\in\mathbf{C}^{n}(G,\mathcal{A}%
L)$, then the mapping $\varphi_{\omega}:\mathcal{A}^{n}G\times
D^{n}\rightarrow L$ defined by 
\begin{equation*}
\varphi_{\omega}(\gamma,d_{1},...,d_{n})=\omega(\gamma)(d_{1}...d_{n}) 
\end{equation*}
abides by the following conditions:

\begin{enumerate}
\item  We have 
\begin{equation*}
\varphi_{\omega}(a\underset{i}{\cdot}\gamma,d_{1},...,d_{n})=\varphi_{\omega
}(\gamma,ad_{1},d_{2},...,d_{n})\text{ \ \ \ \ }(1\leq i\leq n) 
\end{equation*}
for any $a\in\mathbb{R}$.

\item  We have 
\begin{align*}
&
\varphi_{\omega}(%
\gamma,d_{1},...,d_{i-1},ad_{i},d_{i+1},...d_{j-1},d_{j},d_{j+1},...,d_{n})
\\
&
=\varphi_{\omega}(%
\gamma,d_{1},...,d_{i-1},d_{i},d_{i+1},...d_{j-1},ad_{j},d_{j+1},...,d_{n})%
\text{ \ \ \ \ }(1\leq i<j\leq n)
\end{align*}
for any $a\in\mathbb{R}$.

\item  We have 
\begin{equation*}
\varphi_{\omega}(\gamma\circ D^{\sigma},d_{1},...,d_{n})=\varphi_{\omega
}(\gamma,d_{1},...,d_{n})^{\mathrm{sign}(\sigma)}
\end{equation*}
for any permutation $\sigma$ of $\{1,...,n\}$.
\end{enumerate}

Conversely we have

\begin{theorem}
\label{t1.3.2}If $\varphi:\mathcal{A}^{n}G\times D^{n}\rightarrow L$
satisfies the above three conditions, then there exists a unique $%
\omega_{\varphi}\in\mathbf{C}^{n}(G,\mathcal{A}L)$ such that 
\begin{equation*}
\varphi(\gamma,d_{1},...,d_{n})=\omega_{\varphi}(\gamma)(d_{1}...d_{n}) 
\end{equation*}
for any $\gamma\in\mathcal{A}^{n}G$ and any $(d_{1},...,d_{n})\in D^{n}$.
\end{theorem}

\begin{proof}
By the same token as in the proof of Proposition 2 of \S4.2 of Lavendhomme 
\cite{l1} except for the following quasi-colimit diagram of small objects in
place of the corresponding one given there.
\end{proof}

\begin{lemma}
The diagram 
\begin{equation*}
D^{n+1} 
\begin{array}{c}
\overset{\tau_{1}}{\rightarrow} \\ 
\vdots \\ 
\overset{\tau_{n}}{\rightarrow}
\end{array}
D^{n}\overset{\mathfrak{m}}{\rightarrow}D 
\end{equation*}
is a quasi-colimit diagram of small objects, where 
\begin{equation*}
\tau_{i}(d_{0},d_{1},...,d_{n})=(d_{1},...,d_{i-1},d_{0}d_{i},d_{i+1},...d_{n})%
\text{ \ \ \ \ \ }(1\leq i\leq n) 
\end{equation*}
for any $(d_{0},d_{1},...,d_{n})\in D^{n+1}$, and 
\begin{equation*}
\mathfrak{m(}d_{1},...,d_{n})=d_{1}...d_{n}
\end{equation*}
for any $(d_{1},...,d_{n})\in D^{n}$.
\end{lemma}

\subsection{Nishimura Algebroids}

The notion of Nishimura algebroid was introduced by H. Nishimura \cite{n4}
so as to replace the familiar notion of Lie algebroid. In this paper we do
not need its full power, so that our present explanation on Nishimura
algebroid is apparently ad hoc. In this paper, by a \textit{Nishimura
algbroid over a given microlinear space }$M$\textit{,} we mean a pair $%
\mathcal{A=}(\mathcal{A}^{1},\mathcal{A}^{2})$ of microlinear spaces
together with mappings $\pi _{1}:\mathcal{A}^{1}\rightarrow M$, $\pi_{2}:%
\mathcal{A}^{2}\rightarrow M$, $\mathbf{a}_{1}:\mathcal{A}^{1}\rightarrow
M^{D}$, $\mathbf{a}_{2}:\mathcal{A}^{2}\rightarrow M^{D^{2}}$, $\mathbf{s}%
_{i}:\mathcal{A}^{1}\rightarrow\mathcal{A}^{2}$ and $\mathbf{d}_{i}:\mathcal{%
A}^{2}\rightarrow\mathcal{A}^{1}$ ($i=1,2$) such that

\begin{enumerate}
\item  The mappings $\mathbf{s}_{i}:\mathcal{A}^{1}\rightarrow\mathcal{A}^{2}
$ and $\mathbf{d}_{i}:\mathcal{A}^{2}\rightarrow\mathcal{A}^{1}$ ($i=1,2 $)
abide by the so-called simplicial identities. We have $\pi_{2}\circ \mathbf{s%
}_{i}=\pi_{1}$ and $\pi_{1}\circ\mathbf{d}_{i}=\pi_{2}$, while we have $%
\mathbf{d}_{i}\circ\mathbf{a}_{2}=\mathbf{a}_{1}\circ\mathbf{d}_{i}$ and $%
\mathbf{s}_{i}\circ\mathbf{a}_{1}=\mathbf{a}_{2}\circ\mathbf{s}_{i}$, where $%
\mathbf{d}_{1}:M^{D^{2}}\rightarrow M^{D}$ and $\mathbf{s}%
_{i}:M^{D}\rightarrow M^{D^{2}}$ on the left hands of both identities denote
the mappings $\mathbf{d}_{i}:\mathcal{A}^{2}G\rightarrow\mathcal{A}^{1}G$ $%
\mathbf{s}_{i}:\mathcal{A}^{1}G\rightarrow\mathcal{A}^{2}G$ depicted in
Subsection \ref{p1} in case of the pair groupoid $G=M\times M$.

\item  The mapping $\pi_{1}:\mathcal{A}^{1}\rightarrow M$ is a vector bundle.

\item  The square 
\begin{equation*}
\begin{array}{ccc}
\mathcal{A}^{2} & 
\begin{array}{c}
\mathbf{d}_{1} \\ 
\rightarrow
\end{array}
& \mathcal{A}^{1} \\ 
\begin{array}{cc}
\mathbf{d}_{2} & \downarrow
\end{array}
&  & 
\begin{array}{cc}
\downarrow & \pi_{1}
\end{array}
\\ 
\mathcal{A}^{1} & 
\begin{array}{c}
\mathbf{\rightarrow} \\ 
\pi_{1}
\end{array}
& M
\end{array}
\end{equation*}
is a double vector bundle over $M$ in the sense of Mackenzie \cite{ma1},
\S9.1.

\item  A mapping from $\{(\zeta ,x)\in (\mathcal{A}^{1})^{D}\times \mathcal{A%
}^{1}\mid (\pi _{1})^{D}(\zeta )=\mathbf{a}_{1}(x)\}$ to $\mathcal{A}^{2}$,
to be denoted by $\ast $, is defined, where $(\pi _{1})^{D}(\zeta )\in M^{D}$
assigns $\pi _{1}(\zeta (d))$ to each $d\in D$. Any element of $\mathcal{A}%
^{2}$ is uniquely expressed as $\zeta \ast x$. We require that
\begin{equation*}
\mathbf{a}_{2}(\zeta \ast x)=(\mathbf{a}_{1})^{D}(\zeta )
\end{equation*}
where $(\mathbf{a}_{1})^{D}(\zeta )$ assigns $\mathbf{a}_{1}(\zeta
(d_{1}))(d_{2})$ to each $(d_{1},d_{2})\in D^{2}$.
\end{enumerate}

Every groupoid $G$ over $M$ naturally gives rise to a Nishimura algebroid $%
\mathcal{A}G=(\mathcal{A}^{1}G,\mathcal{A}^{2}G)$, where, for any $%
(\zeta,t)\in(\mathcal{A}^{1})^{D}\times\mathcal{A}^{1}$ with $%
(\pi_{1})^{D}(\zeta)=\mathbf{a}_{1}(t)$, $\zeta\ast x$ is defined to be 
\begin{equation*}
(\zeta\ast t)(d_{1},d_{2})=\zeta(d_{1})_{d_{2}}t_{d_{1}}
\end{equation*}
for any $(d_{1},d_{2})\in D^{2}$.

The notion of a \textit{homomorphism} between Nishimura algebroids $\mathcal{%
A=}(\mathcal{A}^{1},\mathcal{A}^{2})$ and $\mathcal{B=}(\mathcal{B}^{1},%
\mathcal{B}^{2})$ over the same base space $M$\ is a pair $%
\varphi=(\varphi_{1},\varphi_{2})$ of mappings $\varphi_{1}:\mathcal{A}%
^{1}\rightarrow\mathcal{B}^{1}$ and $\varphi_{2}:\mathcal{A}^{2}\rightarrow%
\mathcal{B}^{2} $ preserving the structures depicted in the above four
conditions. In particular, we require that 
\begin{equation*}
\varphi_{2}(\zeta\ast x)=(\varphi_{1})^{D}(\zeta)\ast\varphi_{1}(x) 
\end{equation*}
for any $(\zeta,x)\in(\mathcal{A}^{1})^{D}\times\mathcal{A}^{1}$ with $%
(\pi_{1})^{D}(\zeta)=\mathbf{a}_{1}(x)$. Given a Nishimura algebroid $%
\mathcal{A}$ and a vector bundle $E$ over the same base space $M$, a \textit{%
representation} of $\mathcal{A}$ in $E$ is a homomorphism of $\mathcal{A}$
into the Nishimura algebroid $(\mathcal{A}^{1}(\Phi _{\mathrm{Lin}}(E)),%
\mathcal{A}^{2}(\Phi_{\mathrm{Lin}}(E)))$, where $\Phi_{\mathrm{Lin}}(E)$ is
the linear frame groupoid of $E$.

\section{\label{s2}The Additive Complex}

Let $G$ be a groupoid and $E$ a vector bundle over the same base space $E$%
.Let $\rho:\mathcal{A}G\rightarrow\mathcal{A(}\Phi_{\mathrm{Lin}}(E))$ be a
representation of the Nishimura algebroid $\mathcal{A}G$ of $G$ on $E$.
These entities shall be fixed throughout this section. Given $\gamma\in 
\mathcal{A}^{n+1}G$ and $e\in D$, we define $\gamma_{e}^{i}\in\mathcal{A}%
^{n}G\mathcal{\ }$($1\leq i\leq n+1$) to be

\begin{equation*}
\gamma_{e}^{i}(d_{1},...,d_{n})=\gamma(d_{1},...,d_{i-1},e,d_{i},...,d_{n})%
\gamma(0,...,0,\underset{i}{e},0,...,0)^{-1}
\end{equation*}
for any $(d_{1},...,d_{n})\in D^{n}$. Similarly, given $\gamma\in \mathcal{A}%
^{n+1}G$, we define $\gamma_{i}\in\mathcal{A}G$ ($1\leq i\leq n+1$) to be 
\begin{equation*}
\gamma_{i}(d)=\gamma(0,...,0,\underset{i}{d},0,...,0) 
\end{equation*}
for any $d\in D$.

\begin{theorem}
\label{t4.1}Given $\omega\in\mathbf{C}^{n}(G,E)$, there exists a unique $%
\mathbf{d}_{+}\omega\in\mathbf{C}^{n+1}(G,E)$ such that 
\begin{align}
& d_{1}...d_{n+1}\mathbf{d}_{+}\omega(\gamma)  \notag \\
& =\tsum _{i=1}^{n+1}(-1)^{i}d_{1}...\widehat{d_{i}}...d_{n+1}\{\omega(%
\gamma_{0}^{i})-(\rho
(\gamma_{i})_{d_{i}})^{-1}(\omega(\gamma_{d_{i}}^{i}))\}   \label{2.1}
\end{align}
for any $\gamma\in\mathcal{A}^{n+1}G$ and any $(d_{1},...,d_{n+1})\in
D^{n+1} $.
\end{theorem}

\begin{proof}
By Theorem \ref{t1.3.1}, it suffices to note that the function $\varphi :%
\mathcal{A}^{n+1}G\times D^{n+1}\rightarrow E$ defined by 
\begin{align*}
& \varphi(\gamma,d_{1},...,d_{n+1}) \\
& =\tsum _{i=1}^{n+1}(-1)^{i}d_{1}...\widehat{d_{i}}...d_{n+1}\{\omega(%
\gamma_{0}^{i})-(\rho
(\gamma_{i})_{d_{i}})^{-1}(\omega(\gamma_{d_{i}}^{i}))\}
\end{align*}
for any $\gamma\in\mathcal{A}^{n+1}G$ and any $(d_{1},...,d_{n+1})\in D^{n+1}
$ satisfies the three conditions mentioned therein. We fix a notation in
passing. Let $F_{i}:D\rightarrow E$ be the assignment of $(\rho(\gamma
_{i})_{d_{i}})^{-1}(\omega(\gamma_{d}^{i}))\in E$ to each $d\in D$, for
which we have 
\begin{equation*}
\omega(\gamma_{0}^{i})-(\rho(\gamma_{i})_{d_{i}})^{-1}(\omega(%
\gamma_{d_{i}}^{i}))=-d\mathbf{D}F_{i}
\end{equation*}
where $\mathbf{D}F_{i}$ is the derivative of $F_{i}$ at $0$.

\begin{enumerate}
\item  For the first condition, we have to show that 
\begin{equation*}
\varphi(a\underset{j}{\cdot}\gamma,d_{1},...,d_{n})=a\varphi(\gamma
,d_{1},...,d_{n})\text{ \ \ \ \ }(1\leq j\leq n) 
\end{equation*}
for any $a\in\mathbb{R}$. If $i<j$, then we have 
\begin{align*}
& \omega((a\underset{j}{\cdot}\gamma)_{0}^{i})-(\rho(a\underset{j}{\cdot }%
\gamma)_{d_{i}})^{-1}(\omega((a\underset{j}{\cdot}\gamma)_{d_{i}}^{i})) \\
& =\omega((a\underset{j-1}{\cdot}(\gamma_{0}^{i}))-(\rho(%
\gamma_{i})_{d_{i}})^{-1}(\omega((a\underset{j-1}{\cdot}(\gamma_{d_{i}}^{i}))
\\
& =a\{\omega(\gamma_{0}^{i})-(\rho(\gamma_{i})_{d_{i}})^{-1}(\omega
(\gamma_{d_{i}}^{i}))\}
\end{align*}
If $j<i$, then we have 
\begin{align*}
& \omega((a\underset{j}{\cdot}\gamma)_{0}^{i})-(\rho((a\underset{j}{\cdot }%
\gamma)_{i})_{_{d_{i}}})^{-1}(\omega((a\underset{j}{\cdot}%
\gamma)_{d_{i}}^{i})) \\
& =\omega((a\underset{j}{\cdot}(\gamma_{0}^{i}))-(\rho(%
\gamma_{i})_{d_{i}})^{-1}(\omega((a\underset{j}{\cdot}(\gamma_{d_{i}}^{i}))
\\
& =a\{\omega(\gamma_{0}^{i})-(\rho(\gamma_{i})_{d_{i}})^{-1}(\omega
(\gamma_{d_{i}}^{i}))\}
\end{align*}
Finally we consider the case of $j=i$, in which we have 
\begin{align*}
& \omega((a\underset{i}{\cdot}\gamma)_{0}^{i})-(\rho((a\underset{i}{\cdot }%
\gamma)_{i})_{_{d_{i}}})^{-1}(\omega((a\underset{i}{\cdot}%
\gamma)_{d_{i}}^{i})) \\
& =\omega(\gamma_{0}^{i})-(\rho(\gamma_{i})_{ad_{i}})^{-1}(\omega
(\gamma_{ad_{i}}^{i})) \\
& =-ad_{i}\mathbf{D}F_{i} \\
& =a\{\omega(\gamma_{0}^{i})-(\rho(\gamma_{i})_{d_{i}})^{-1}(\omega
(\gamma_{d_{i}}^{i}))\}
\end{align*}

\item  The easy verification of the second condition is left to the reader.

\item  For the third condition, it suffices to note that 
\begin{align*}
& d_{1}...\widehat{d_{i}}...d_{n+1}\{\omega(\gamma_{0}^{i})-(\rho(\gamma
_{i})_{d_{i}})^{-1}(\omega(\gamma_{d_{i}}^{i}))\} \\
& =-d_{1}...d_{n+1}\mathbf{D}F_{i} \\
& =d_{1}...\widehat{d_{j}}...d_{n+1}\{\omega(\gamma_{0}^{i})-(\rho(\gamma
_{i})_{d_{j}})^{-1}(\omega(\gamma_{d_{j}}^{i}))\}
\end{align*}
for any $i\neq j$.
\end{enumerate}
\end{proof}

Now we would like to show that $\mathbf{d}_{+}^{2}=0$, for which we need
three lemmas. Let $\gamma\in\mathcal{A}^{n+2}G$ throughout the following
three lemmas.

\begin{lemma}
\label{t2.1}For $1\leq j\leq i\leq n+1$, we have 
\begin{equation*}
(\gamma_{0}^{j})_{i}=\gamma_{i+1}
\end{equation*}
For $1\leq i<j\leq n+2$ we have 
\begin{equation*}
(\gamma_{0}^{j})_{i}=\gamma_{i}
\end{equation*}
\end{lemma}

\begin{proof}
Obvious.
\end{proof}

\begin{lemma}
\label{t4.2}For $1\leq j<i\leq n+2$ and $e,e^{\prime}\in D$, we have 
\begin{equation*}
(\gamma_{e}^{i})_{e^{\prime}}^{j}=(\gamma_{e^{\prime}}^{j})_{e}^{i-1}
\end{equation*}
\end{lemma}

\begin{proof}
It is easy to see that both $(\gamma_{e}^{i})_{e^{\prime}}^{j}$ and $%
(\gamma_{e^{\prime}}^{j})_{e}^{i-1}$ are the same mapping as follows: 
\begin{align*}
(d_{1},...,d_{n}) & \in D^{n} \\
&
\longmapsto\gamma(d_{1},...,d_{j-1},e^{%
\prime},d_{j},...,d_{i-1},e,d_{i},...,d_{n}) \\
& \gamma(0,...,0,\underset{j}{e^{\prime}},0,...,0,\underset{i}{e}%
,0,...,0)^{-1}
\end{align*}
\end{proof}

\begin{lemma}
\label{t4.3}For $1\leq j<i\leq n+2$ and $d_{i},d_{j}\in D$, we have 
\begin{equation*}
(\rho(\gamma_{i})_{d_{i}})^{-1}\circ(\rho((%
\gamma_{d_{i}}^{i})_{j})_{d_{j}})^{-1}=(\rho(\gamma_{j})_{d_{j}})^{-1}\circ(%
\rho((\gamma_{d_{j}}^{j})_{i-1})_{d_{i}})^{-1}
\end{equation*}
\end{lemma}

\begin{proof}
It is easy to see that both $(\rho (\gamma _{i})_{d_{i}})^{-1}\circ (\rho
((\gamma _{d_{i}}^{i})_{j})_{d_{j}})^{-1}$ and $(\rho (\gamma
_{j})_{d_{j}})^{-1}\circ (\rho ((\gamma _{d_{j}}^{j})_{i-1})_{d_{i}})^{-1}$
are equal to
\begin{equation*}
(\rho (\gamma (0,...,0,\underset{}{\underset{j}{\cdot _{1}},0,...,0,%
\underset{i}{\cdot _{2}}},0,...,0))_{d_{j},d_{i}})^{-1}
\end{equation*}
To see this, we should note that 
\begin{align*}
& \rho ((\gamma _{d_{i}}^{i})_{j})_{d_{j}}\circ \rho (\gamma _{i})_{d_{i}} \\
& =\{(d\in D\mapsto \rho ((\gamma _{d}^{i})_{j}))\ast \rho (\gamma
_{i})\}_{d_{i},d_{j}} \\
& =\rho ((d\in D\mapsto (\gamma _{d}^{i})_{j})\ast \gamma _{i})_{d_{i},d_{j}}
\\
& =\rho (\gamma (0,...,0,\underset{}{\underset{j}{\cdot _{1}},0,...,0,%
\underset{i}{\cdot _{2}}},0,...,0))_{d_{j},d_{i}}
\end{align*}
and 
\begin{align*}
& \rho ((\gamma _{d_{j}}^{j})_{i-1})_{d_{i}}\circ \rho (\gamma _{j})_{d_{j}}
\\
& =\{(d\in D\mapsto \rho ((\gamma _{d}^{j})_{i-1}))\ast \rho (\gamma
_{j})\}_{d_{j},d_{i}} \\
& =\rho ((d\in D\mapsto (\gamma _{d}^{j})_{i-1})\ast \gamma
_{j})_{d_{j},d_{i}} \\
& =\rho (\gamma (0,...,0,\underset{}{\underset{j}{\cdot _{1}},0,...,0,%
\underset{i}{\cdot _{2}}},0,...,0))_{d_{j},d_{i}}
\end{align*}
\end{proof}

\begin{theorem}
\label{t4.4}We have 
\begin{equation*}
\mathbf{d}_{+}^{2}=0\text{.}
\end{equation*}
In other words, the composition 
\begin{equation*}
\mathbf{C}^{n}(G,E)\overset{\mathbf{d}_{+}}{\rightarrow}\mathbf{C}^{n+1}(G,E)%
\overset{\mathbf{d}_{+}}{\rightarrow}\mathbf{C}^{n+2}(G,E) 
\end{equation*}
vanishes.
\end{theorem}

\begin{proof}
Let $\gamma\in\mathcal{A}^{n+2}G$, $d_{1},...,d_{n+2}\in D$ and $\omega \in%
\mathbf{C}^{n}(G,E)$. We have 
\begin{align*}
& d_{1}...d_{n+2}\mathbf{d}_{+}^{2}\omega(\gamma) \\
& =\tsum _{i=1}^{n+2}(-1)^{i}d_{1}...\widehat{d_{i}}...d_{n+2}\mathbf{d}%
_{+}\omega(\gamma_{0}^{i})- \\
& \tsum _{i=1}^{n+2}(-1)^{i}d_{1}...\widehat{d_{i}}...d_{n+2}(\rho(%
\gamma_{i})_{d_{i}})^{-1}(\mathbf{d}_{+}\omega(\gamma_{d_{i}}^{i})) \\
& =\tsum _{i=1}^{n+2}\tsum _{j=1}^{i-1}(-1)^{i+j}d_{1}...\widehat{d_{j}}...%
\widehat{d_{i}}...d_{n+2}\omega(\left( \gamma_{0}^{i}\right) _{0}^{j})+ \\
& \tsum _{i=1}^{n+2}\tsum _{j=1}^{i-1}(-1)^{i+j+1}d_{1}...\widehat{d_{j}}...%
\widehat{d_{i}}...d_{n+2}(\rho
((\gamma_{0}^{i})_{j})_{d_{j}})^{-1}(\omega(\left( \gamma_{0}^{i}\right)
_{d_{j}}^{j}))+ \\
& \tsum _{i=1}^{n+2}\tsum _{j=1}^{i-1}(-1)^{i+j+1}d_{1}...\widehat{d_{j}}...%
\widehat{d_{i}}...d_{n+2}(\rho
(\gamma_{i})_{d_{i}})^{-1}(\omega((\gamma_{d_{i}}^{i})_{0}^{j}))+ \\
& \tsum _{i=1}^{n+2}\tsum _{j=1}^{i-1}(-1)^{i+j}d_{1}...\widehat{d_{j}}...%
\widehat{d_{i}}...d_{n+2}((\rho(\gamma
_{i})_{d_{i}})^{-1}\circ(\rho((\gamma_{d_{i}}^{i})_{j})_{d_{j}})^{-1})(%
\omega((\gamma_{d_{i}}^{i})_{d_{j}}^{j}))+ \\
& \tsum _{i=1}^{n+2}\tsum _{j=i}^{n+1}(-1)^{i+j}d_{1}...\widehat{d_{i}}...%
\widehat{d_{j+1}}...d_{n+2}\omega(\left( \gamma_{0}^{i}\right) _{0}^{j})+ \\
& \tsum _{i=1}^{n+2}\tsum _{j=i}^{n+1}(-1)^{i+j+1}d_{1}...\widehat{d_{i}}...%
\widehat{d_{j+1}}...d_{n+2}(\rho
((\gamma_{0}^{i})_{j})_{d_{j}})^{-1}\omega(\left( \gamma_{0}^{i}\right)
_{d_{j}}^{j}))+ \\
& \tsum _{i=1}^{n+2}\tsum _{j=i}^{n+1}(-1)^{i+j+1}d_{1}...\widehat{d_{i}}...%
\widehat{d_{j+1}}...d_{n+2}(\rho
(\gamma_{i})_{d_{i}})^{-1}(\omega((\gamma_{d_{i}}^{i})_{0}^{j}))+ \\
& \tsum _{i=1}^{n+2}\tsum _{j=i}^{n+1}(-1)^{i+j}d_{1}...\widehat{d_{i}}...%
\widehat{d_{j+1}}...d_{n+2}((\rho
(\gamma_{i})_{d_{i}})^{-1}\circ(\rho((%
\gamma_{d_{i}}^{i})_{j})_{d_{j}})^{-1})(\omega((%
\gamma_{d_{i}}^{i})_{d_{j}}^{j})) \\
& =0
\end{align*}
The final derivation of total vanishment comes from the cancellation of the
first double summantion and the fifth one in the previous development by
Lemma \ref{t4.2}, that of the second double summation and the seventh one in
the previous development by Lemmas \ref{t4.1} and \ref{t4.2}, that of the
third double summation and the sixth one in the previous development by
Lemmas \ref{t4.1} and \ref{t4.2}, and finally that of the fourth double
summation and the eighth one in the previous development by Lemmas \ref{t4.2}
and \ref{t4.3}.
\end{proof}

The operator $\mathbf{d}_{+}$ is called the \textit{additive coboundary
operator with respect to the representation} $\rho$.

\section{\label{s3}The Multiplicative Complex}

Let $G$ be a groupoid and $L$ a group bundle over the same base space $E$.
Let $\rho:\mathcal{A}G\rightarrow\mathcal{A(}\Phi_{\mathrm{Lin}}(\mathcal{A}%
^{1}L))$ be a representation of the Nishimura algebroid $\mathcal{A}G$ of $G$
on the vector bundle $\mathcal{A}^{1}L$. These entities shall be fixed
throughout the rest of this paper. The preceding section tells us that we
have the additive coboundary operator 
\begin{equation*}
\mathbf{d}_{+}:\mathbf{C}^{n}(G,\mathcal{A}^{1}L)\overset{}{\rightarrow }%
\mathbf{C}^{n+1}(G,\mathcal{A}^{1}L) 
\end{equation*}
with respect to the representation $\rho:\mathcal{A}G\rightarrow \mathcal{A(}%
\Phi_{\mathrm{Lin}}(\mathcal{A}^{1}L))$. Now we are going to define another
coboundary operator 
\begin{equation*}
\mathbf{d}_{\times}:\mathbf{C}^{n}(G,\mathcal{A}^{1}L)\overset{}{\rightarrow 
}\mathbf{C}^{n+1}(G,\mathcal{A}^{1}L) 
\end{equation*}
to be called the \textit{multiplicative coboundary operator with respect to }%
$\rho$. Given $\omega\in\mathbf{C}^{n}(G,\mathcal{A}^{1}L)$, $\mathbf{d}%
_{\times}\omega\in\mathbf{C}^{n+1}(G,\mathcal{A}^{1}L)$ is expected to be
defined in such a way that 
\begin{align}
& ((\mathbf{d}_{\times}\omega)(\gamma))_{d_{1}...d_{n+1}}  \notag \\
& =\prod_{i=1}^{n+1}\{(\omega(\gamma_{0}^{i}))_{d_{1}...\widehat{d_{i}}%
...d_{n+1}}(\rho(\gamma_{i})_{d_{i}})^{-1}((\omega(%
\gamma_{d_{i}}^{i}))_{-d_{1}...\widehat{d_{i}}...d_{n+1}})\}^{(-1)^{i}} 
\label{3.1}
\end{align}
for any $\gamma\in\mathcal{A}^{n+1}G$ and any $(d_{1},...,d_{n+1})\in D^{n+1}
$. To show its existence and uniqueness, we need a simple lemma.

\begin{lemma}
\label{t3.1}Let $F_{i}:D\rightarrow\mathcal{A}^{1}L$ be the assignment in
the proof of Theorem \ref{t4.1} with $\mathcal{A}^{1}L$ in place of $E$.
Then we have 
\begin{align*}
& (\omega(\gamma_{0}^{i}))_{-d}(\rho(\gamma_{i})_{e})^{-1}((\omega(\gamma
_{e}^{i}))_{d}) \\
& =(e\mathbf{D}F_{i})_{d} \\
& =(\mathbf{D}F_{i})_{de}
\end{align*}
\end{lemma}

\begin{proof}
By dint of Proposition \ref{t1.1.3}, the statement is merely a reformulation
of 
\begin{equation*}
F_{i}(e)=F_{i}(0)+e\mathbf{D}F_{i}
\end{equation*}
\end{proof}

\begin{corollary}
\label{t3.2}The order of multiplication of $n+1$ factors in (\ref{3.1}) does
not matter.
\end{corollary}

\begin{proof}
By dint of Proposition \ref{t1.1.3}, it suffices to note that 
\begin{align*}
& \{(\omega(\gamma_{0}^{i}))_{d_{1}...\widehat{d_{i}}...d_{n+1}}(\rho
(\gamma_{i})_{d_{i}})^{-1}((\omega(\gamma_{d_{i}}^{i}))_{-d_{1}...\widehat
{d_{i}}...d_{n+1}})\}^{(-1)^{i}} \\
& =(\mathbf{D}F_{i})_{(-1)^{i+1}d_{1}...d_{n+1}}
\end{align*}
\end{proof}

\begin{theorem}
\label{t3.3}For any $\omega\in\mathbf{C}^{n}(G,\mathcal{A}^{1}L)$, there
exists a unique $\mathbf{d}_{\times}\omega\in\mathbf{C}^{n+1}(G,\mathcal{A}%
^{1}L)$ abiding by the condition (\ref{3.1}).
\end{theorem}

\begin{proof}
By Theorem \ref{t1.3.2}, it suffices to note that the function $\varphi :%
\mathcal{A}^{n+1}G\times D^{n+1}\rightarrow L$ defined by 
\begin{align*}
& \varphi(\gamma,d_{1},...,d_{n+1}) \\
& =\prod_{i=1}^{n+1}\{(\omega(\gamma_{0}^{i}))_{d_{1}...\widehat{d_{i}}%
...d_{n+1}}(\rho(\gamma_{i})_{d_{i}})^{-1}((\omega(%
\gamma_{d_{i}}^{i}))_{-d_{1}...\widehat{d_{i}}...d_{n+1}})\}^{(-1)^{i}}
\end{align*}
for any $\gamma\in\mathcal{A}^{n+1}G$ and any $(d_{1},...,d_{n+1})\in D^{n+1}
$ satisfies the three conditions mentioned therein. The proof can be carried
out as in the proof of Theorem \ref{t2.1} by dint of Lemma \ref{t3.1}.
\end{proof}

It is to be shown in the succeeding section that $\mathbf{d}_{\times}^{2}=0$.

\section{The Coincidence of the Two Complexes}

Now we come to show the main result of this paper.

\begin{theorem}
\label{t5.1}(The Coincidence Theorem). The additive coboundary operator $%
\mathbf{d}_{+}:\mathbf{C}^{n}(G,\mathcal{A}^{1}L)\rightarrow\mathbf{C}%
^{n+1}(G,\mathcal{A}^{1}L)$\ with respect to the representation $\rho :%
\mathcal{A}G\rightarrow\mathcal{A(}\Phi_{\mathrm{Lin}}(\mathcal{A}^{1}L))$
and the multiplicative coboundary operator $\mathbf{d}_{\times}:\mathbf{C}%
^{n}(G,\mathcal{A}^{1}L)\rightarrow\mathbf{C}^{n+1}(G,\mathcal{A}^{1}L)$\
with respect to\textit{\ the representation }$\rho:\mathcal{A}G\rightarrow 
\mathcal{A(}\Phi_{\mathrm{Lin}}(\mathcal{A}^{1}L))$ coincide for any natural
number $n$.
\end{theorem}

\begin{proof}
Let $\gamma\in\mathcal{A}^{n+1}G$ and $e_{1},...,e_{n+1},d_{1},...,d_{n+1}%
\in D$. We note that 
\begin{align*}
& (e_{1}...e_{n+1}(\mathbf{d}_{\times}\omega(\gamma)))_{d_{1}...d_{n+1}} \\
& =(\mathbf{d}_{\times}\omega(\gamma))_{d_{1}e_{1}...d_{n+1}e_{n+1}} \\
& =\prod_{i=1}^{n+1}((\omega(\gamma_{0}^{i}))_{d_{1}e_{1}...\widehat
{d_{i}e_{i}}...d_{n+1}e_{n+1}}(\rho(\gamma_{i})_{d_{i}})^{-1}((\omega
(\gamma_{d_{i}e_{i}}^{i}))_{-d_{1}e_{1}...\widehat{d_{i}e_{i}}%
...d_{n+1}e_{n+1}}))^{(-1)^{i}} \\
& =\prod_{i=1}^{n+1}((d_{i}e_{i}\mathbf{D}F_{i})_{-d_{1}e_{1}...\widehat
{d_{i}e_{i}}...d_{n+1}e_{n+1}})^{(-1)^{i}} \\
& =\prod_{i=1}^{n+1}((\mathbf{D}%
F_{i})_{-d_{1}e_{1}...d_{n+1}e_{n+1}})^{(-1)^{i}} \\
& =\prod_{i=1}^{n+1}((e_{1}...e_{n+1}\mathbf{D}%
F_{i})_{-_{d_{1}...d_{n+1}}})^{(-1)^{i}} \\
& =\{\tsum _{i=1}^{n+1}(-1)^{i}(e_{1}...\widehat{e_{i}}...e_{n+1}\omega(%
\gamma_{0}^{i})-e_{1}...\widehat{e_{i}}...e_{n+1}(\rho(%
\gamma_{i})_{e_{i}})^{-1}(\omega (\gamma_{e_{i}}^{i})))\}_{d_{1}...d_{n+1}}
\\
& =(e_{1}...e_{n+1}(\mathbf{d}_{+}\omega(\gamma)))_{d_{1}...d_{n+1}}
\end{align*}
Since $e_{1},...,e_{n+1},d_{1},...,d_{n+1}\in D$ were arbitrary, we can
conclude that $\mathbf{d}_{+}\omega(\gamma)=\mathbf{d}_{\times}\omega(\gamma)
$.
\end{proof}

\begin{corollary}
\label{t5.2}We have 
\begin{equation*}
\mathbf{d}_{\times}^{2}=0\text{.}
\end{equation*}
In other words, the composition 
\begin{equation*}
\mathbf{C}^{n}(G,\mathcal{A}^{1}L)\overset{\mathbf{d}_{\times}}{\rightarrow }%
\mathbf{C}^{n+1}(G,\mathcal{A}^{1}L)\overset{\mathbf{d}_{\times}}{\rightarrow%
}\mathbf{C}^{n+2}(G,\mathcal{A}^{1}L) 
\end{equation*}
vanishes.
\end{corollary}

\begin{proof}
This follows simply from Theorems \ref{t4.4} and \ref{t5.1}.
\end{proof}

\section{The Comparison between the Multiplicative Coboundary Operator and
Kock's Contour Derivative from Dimension $1$ to Dimension $2$}

Kock's contour derivative $\mathbf{d}_{\mathrm{\circlearrowleft}}:\mathbf{C}%
^{0}(G,\mathcal{A}^{1}L)\overset{}{\rightarrow}\mathbf{C}^{1}(G,\mathcal{A}%
^{1}L)$ should unquestionably agree with the multiplicative coboundary
operator $\mathbf{d}_{\times}:\mathbf{C}^{0}(G,\mathcal{A}^{1}L)\overset{}{%
\rightarrow}\mathbf{C}^{1}(G,\mathcal{A}^{1}L)$, for which there is nothing
to discuss. Kock's succeeding contour derivative $\mathbf{d}_{\mathrm{%
\circlearrowleft}}:\mathbf{C}^{1}(G,\mathcal{A}^{1}L)\overset{}{\rightarrow}%
\mathbf{C}^{2}(G,\mathcal{A}^{1}L)$ should diverge from $\mathbf{d}_{\times}:%
\mathbf{C}^{1}(G,\mathcal{A}^{1}L)\overset{}{\rightarrow}\mathbf{C}^{2}(G,%
\mathcal{A}^{1}L)$, and it is our main concern here how $\mathbf{d}_{\mathrm{%
\circlearrowleft}}:\mathbf{C}^{1}(G,\mathcal{A}^{1}L)\overset{}{\rightarrow}%
\mathbf{C}^{2}(G,\mathcal{A}^{1}L)$ diverges exactly from $\mathbf{d}%
_{\times}:\mathbf{C}^{1}(G,\mathcal{A}^{1}L)\overset{}{\rightarrow}\mathbf{C}%
^{2}(G,\mathcal{A}^{1}L)$.

Given $\gamma\in\mathcal{A}^{2}G$, we define $\mathbf{d}_{\mathrm{%
\circlearrowleft}}\omega(\gamma)\in\mathcal{A}^{1}L$ to be 
\begin{align}
& (\mathbf{d}_{\mathrm{\circlearrowleft}}\omega(\gamma))_{d_{1}d_{2}}  \notag
\\
&
=\omega(\gamma_{0}^{1})_{-d_{2}}\{(\rho(\gamma_{2})_{d_{2}})^{-1}(\omega(%
\gamma_{d_{2}}^{2}))\}_{-d_{1}}\{(\rho(\gamma_{1})_{d_{1}})^{-1}(\omega(%
\gamma_{d_{1}}^{1}))\}_{d_{2}})\omega(\gamma_{0}^{2})_{d_{1}}   \label{6.1}
\end{align}
for $d_{1},d_{2}\in D$. Obviously we have to show that

\begin{proposition}
For any $\omega\in\mathbf{C}^{1}(G,\mathcal{A}^{1}L)$ and any $\gamma \in%
\mathcal{A}^{2}G$, there exists a unique $\mathbf{d}_{\mathrm{%
\circlearrowleft}}\omega(\gamma)\in\mathcal{A}^{1}L$ abiding by the
condition (\ref{6.1}).
\end{proposition}

\begin{proof}
It suffices to note that if $d_{1}=0$ or $d_{2}=0$, then the right-hand side
of (\ref{6.1}) is equal to the identity at $\gamma(0,0)$, which is easy to
see.
\end{proof}

Now we are ready to establish the main result of Kock \cite{k3} in our
general context.

\begin{theorem}
Given $\omega\in\mathbf{C}^{1}(G,\mathcal{A}^{1}L)$ and $\gamma\in \mathcal{A%
}^{2}G$, we have 
\begin{equation*}
\mathbf{d}_{\mathrm{\circlearrowleft}}\omega(\gamma)=\mathbf{d}%
_{\times}\omega(\gamma)+[\omega(\gamma_{0}^{2}),\omega(\gamma_{0}^{1})] 
\end{equation*}
\end{theorem}

\begin{proof}
Let $d_{1},d_{2}\in D$. Then we have 
\begin{align*}
& \{\mathbf{d}_{\mathrm{\circlearrowleft}}\omega(\gamma)-\mathbf{d}_{\times
}\omega(\gamma)\}_{d_{1}d_{2}} \\
& =\{\mathbf{d}_{\times}\omega(\gamma)\}_{-d_{1}d_{2}}\{\mathbf{d}_{\mathrm{%
\circlearrowleft}}\omega(\gamma)\}_{d_{1}d_{2}} \\
&
=[\omega(\gamma_{0}^{2})_{-d_{1}}\{(\rho(\gamma_{2})_{d_{2}})^{-1}(\omega(%
\gamma_{d_{2}}^{2}))\}_{d_{1}}\{(\rho(\gamma_{1})_{d_{1}})^{-1}(\omega(%
\gamma_{d_{1}}^{1}))\}_{-d_{2}}\omega(\gamma_{0}^{1})_{d_{2}}] \\
&
[\omega(\gamma_{0}^{1})_{-d_{2}}\{(\rho(\gamma_{2})_{d_{2}})^{-1}(\omega(%
\gamma_{d_{2}}^{2}))\}_{-d_{1}}\{(\rho(\gamma_{1})_{d_{1}})^{-1}(\omega(%
\gamma_{d_{1}}^{1}))\}_{d_{2}}\omega(\gamma_{0}^{2})_{d_{1}}] \\
&
=\omega(\gamma_{0}^{2})_{-d_{1}}\{(\rho(\gamma_{2})_{d_{2}})^{-1}(\omega(%
\gamma_{d_{2}}^{2}))\}_{d_{1}}\{(\rho(\gamma_{1})_{d_{1}})^{-1}(\omega(%
\gamma_{d_{1}}^{1}))\}_{-d_{2}} \\
&
\{(\rho(\gamma_{2})_{d_{2}})^{-1}(\omega(\gamma_{d_{2}}^{2}))\}_{-d_{1}}\{(%
\rho(\gamma_{1})_{d_{1}})^{-1}(\omega(\gamma_{d_{1}}^{1}))\}_{d_{2}}\omega(%
\gamma_{0}^{2})_{d_{1}}
\end{align*}
Since we have 
\begin{align*}
\omega(\gamma_{0}^{2})_{-d_{1}}\{(\rho(\gamma_{2})_{d_{2}})^{-1}(\omega
(\gamma_{d_{2}}^{2}))\}_{d_{1}} & =(\mathbf{D}F_{2})_{d_{1}d_{2}} \\
\{(\rho(\gamma_{1})_{d_{1}})^{-1}(\omega(\gamma_{d_{1}}^{1}))\}_{-d_{2}} & =(%
\mathbf{D}F_{1})_{-d_{1}d_{2}}\omega(\gamma_{0}^{1})_{-d_{2}} \\
\{(\rho(\gamma_{2})_{d_{2}})^{-1}(\omega(\gamma_{d_{2}}^{2}))\}_{-d_{1}} & =(%
\mathbf{D}F_{2})_{-d_{1}d_{2}}\omega(\gamma_{0}^{1})_{-d_{1}} \\
\{(\rho(\gamma_{1})_{d_{1}})^{-1}(\omega(\gamma_{d_{1}}^{1}))\}_{d_{2}} & =(%
\mathbf{D}F_{1})_{d_{1}d_{2}}\omega(\gamma_{0}^{1})_{d_{2}}
\end{align*}
we can continue our calculation as follows: 
\begin{align*}
& \{\mathbf{d}_{\mathrm{\circlearrowleft}}\omega(\gamma)-\mathbf{d}_{\times
}\omega(\gamma)\}_{d_{1}d_{2}} \\
& =(\mathbf{D}F_{2})_{d_{1}d_{2}}(\mathbf{D}F_{1})_{-d_{1}d_{2}}\omega
(\gamma_{0}^{1})_{-d_{2}}(\mathbf{D}F_{2})_{-d_{1}d_{2}}\omega(%
\gamma_{0}^{1})_{-d_{1}}(\mathbf{D}F_{1})_{d_{1}d_{2}}\omega(%
\gamma_{0}^{1})_{d_{2}}\omega(\gamma_{0}^{2})_{d_{1}}
\end{align*}
It is easy to see that $(\mathbf{D}F_{1})_{d_{1}d_{2}}$, $(\mathbf{D}%
F_{1})_{-d_{1}d_{2}}$, $(\mathbf{D}F_{2})_{d_{1}d_{2}}$ and $(\mathbf{D}%
F_{2})_{-d_{1}d_{2}}$ commute with any term occurring in the above
calculation, so that the calculation itself can step forward as follows with 
$(\mathbf{D}F_{1})_{d_{1}d_{2}}$ and $(\mathbf{D}F_{2})_{d_{1}d_{2}}$
canceling out $(\mathbf{D}F_{1})_{-d_{1}d_{2}}$ and $(\mathbf{D}%
F_{2})_{-d_{1}d_{2}}$ respectively: 
\begin{align*}
& \{\mathbf{d}_{\mathrm{\circlearrowleft}}\omega(\gamma)-\mathbf{d}_{\times
}\omega(\gamma)\}_{d_{1}d_{2}} \\
&
=\omega(\gamma_{0}^{1})_{-d_{2}}\omega(\gamma_{0}^{1})_{-d_{1}}\omega(%
\gamma_{0}^{1})_{d_{2}}\omega(\gamma_{0}^{2})_{d_{1}} \\
& =[\omega(\gamma_{0}^{2}),\omega(\gamma_{0}^{1})]_{d_{1}d_{2}}
\end{align*}
Since $d_{1}\in D$ and $d_{2}\in D$ were arbitrary, the desired conclusion
follows at once.
\end{proof}

\begin{corollary}
For any $\omega\in\mathbf{C}^{1}(G,\mathcal{A}^{1}L)$, we have $\mathbf{d}_{%
\mathrm{\circlearrowleft}}\omega\in\mathbf{C}^{2}(G,\mathcal{A}^{1}L)$.
\end{corollary}

\begin{corollary}
If $\omega\in\mathbf{C}^{1}(G,\mathcal{A}^{1}L)$ is a closed form (i.e., $%
\mathbf{d}_{\times}\omega=0$), then we have 
\begin{equation*}
\mathbf{d}_{\mathrm{\circlearrowleft}}\omega(\gamma)=[\omega(%
\gamma_{0}^{2}),\omega(\gamma_{0}^{1})] 
\end{equation*}
\end{corollary}

\end{document}